\documentclass[12pt,a4paper]{amsart}
\usepackage{amssymb,amsmath}

\headsep=1cm \numberwithin{equation}{section}
\newtheorem{theorem}{Theorem}[section]
\newtheorem{lemma}[theorem]{Lemma}
\newtheorem{proposition}[theorem]{Proposition}
\newtheorem{corollary}[theorem]{Corollary}

\theoremstyle{definition}
\newtheorem{definition}[theorem]{Definition}
\theoremstyle{remark}
\newtheorem{remark}[theorem]{Remark}
\newtheorem{example}[theorem]{Example}

\newcommand{\Hom}{\operatorname{Hom}}

\newcommand{\Ann}{\operatorname{Ann}}

\newcommand{\lo}{\longrightarrow}
\newcommand{\fm}{\frak{m}}

\newcommand{\fa}{\frak{a}}
\newcommand{\fb}{\frak{b}}

\def\mapdown#1{\Big\downarrow\rlap
{$\vcenter{\hbox{$\scriptstyle#1$}}$}}

\begin{document}
\author[Divaani-Aazar, Esmkhani and Tousi ]{Kamran Divaani-Aazar, Mohammad
Ali Esmkhani and Massoud Tousi}
\title[Some criteria of cyclically...]
{Some criteria of cyclically pure injective modules}

\address{K. Divaani-Aazar, Department of Mathematics, Az-Zahra University,
Vanak, Post Code 19834, Tehran, Iran-and-Institute for Studies in
Theoretical Physics and Mathematics, P.O. Box 19395-5746, Tehran,
Iran.} \email{kdivaani@ipm.ir}

\address{M.A. Esmkhani, Department of Mathematics, Shahid Beheshti University,
Tehran, Iran-and-Institute for Studies in Theoretical Physics and
Mathematics, P.O. Box 19395-5746, Tehran, Iran.}

\address{M. Tousi, Department of Mathematics, Shahid Beheshti University,
Tehran, Iran-and-Institute for Studies in Theoretical Physics and
Mathematics, P.O. Box 19395-5746, Tehran, Iran.}


\keywords{Cocyclic modules, cyclic exact sequences, cyclically
pure injective modules, quasi-complete rings, small cofinite
irreducibles.\\
The first author was supported by a grant from IPM (No. 84130213).\\
The third author was supported by a grant from IPM (No.
84130214).}

\begin{abstract} The structure of cyclically pure injective modules
over a commutative ring $R$ is investigated and several
characterizations for them are presented. In particular, we prove
that a module $D$ is cyclically pure injective if and only if $D$
is isomorphic to a direct summand of a module of the form
$\Hom_R(L,E)$ where $L$ is the direct sum of a family of finitely
presented cyclic modules and $E$ is an injective module. Also, we
prove that over a quasi-complete Noetherian ring $(R,\fm)$ an
$R$-module $D$ is cyclically pure injective if and only if there
is a family $\{C_\lambda\}_{\lambda\in \Lambda}$ of cocyclic
modules such that $D$ is isomorphic to a direct summand of
$\Pi_{\lambda\in \Lambda}C_\lambda$. Finally, we show that over a
complete local ring every finitely generated module which has
small cofinite irreducibles is cyclically pure injective.
\end{abstract}

\maketitle

\section{Introduction}

Throughout this paper, let R denote a commutative ring with
identity. All modules are assumed to be left unitary. The notion
of pure injective modules has a substantial role in commutative
algebra and model theory. Even in model theory [{\bf 8}], the
notion of pure injective modules is much more useful than that of
injective modules. Also, there are some nice applications of this
notion in theory of flat covers (see e.g. [{\bf 11}]).

There are several generalizations of the notion of pure injective
modules. One of these generalizations is the notion of cyclically
pure injective modules which has attracted more attention in
recent years. Following his investigations on ``direct summand
conjecture'', M. Hochster [{\bf 6}] studied the structure of
Noetherian rings that are pure in every module in which it is
cyclically pure. He showed that a Noetherian ring $R$ is pure in
every module in which it is cyclically pure if and only if $R$ has
small cofinite irreducibles.  Using the notion of cyclically pure
injective modules, L. Melkersson  [{\bf 7}] provided some
characterizations for a finitely generated module $M$ over a
Noetherain local ring which is pure in every cyclically pure
extension of $M$.  In this paper, our aim is to present some
criterions of cyclically pure injective modules, through a
systematic investigation of their structure.

There are several variants of the notion of purity (see e.g. [{\bf
10}]). More generally, let $\mathcal{S}$ be a class of
$R$-modules. An exact sequence $0\lo A\lo B\lo C\lo 0$ is
$\mathcal{S}$-pure if for all $M\in \mathcal{S}$ the induced
homomorphism $\Hom_R(M,B)\lo \Hom_R(M,C)$ is surjective. An
$R$-module $D$ is said to be $\mathcal{S}$-pure injective if for
any $\mathcal{S}$-pure exact sequence $0\lo A\lo B\lo C\lo 0 ,$
the induced homomorphism $\Hom_R(B,D)\lo \Hom_R(A,D)$ is
surjective. When $\mathcal{S}$ is the class of finitely presented
$R$-modules, $\mathcal{S}$-pure exact sequences and
$\mathcal{S}$-pure injective modules are called pure exact
sequences and pure injective modules, respectively. In this
article, we consider the class $\mathcal{S}$ consisting of all
$R$-modules $M$ for which there are an integer $n\in \mathbb{N}$
and a cyclic submodule $G$ of $R^{n}$ such that $M$ is isomorphic
to $R^{n}/G$. In Section 2, a characterization of cyclically pure
exact sequences is given. Among other things, this
characterization implies that for the above class $\mathcal{S}$,
$\mathcal{S}$-pure exact sequences and $\mathcal{S}$-pure
injective modules coincide with cyclically pure exact sequences
and cyclically pure injective modules, respectively. Also, several
elementary results will be presented in this section, to ease
reading the remainder of the paper.

In Section 3, we present two characterizations of cyclically pure
injective modules. The first one, in particular, asserts that an
$R$-module $D$ is cyclically pure injective if and only if $D$ has
no proper essential cyclically pure extension. Also, it is proved
that an $R$-module $D$ is cyclically pure injective if and only if
$D$ is isomorphic to a direct summand of a module of the form
$\Hom_R(L,E)$ where E is an injective $R$-module and L is the
direct sum of a family of finitely presented cyclic modules.

In Section 4, we show that every $R$-module possesses a unique, up
to isomorphism, cyclically pure injective envelope.

In Section 5, we investigate the question when cocyclic modules
are cyclically pure injective. As a result, we present our last
characterization of pure injective modules. Namely, we prove that
over a quasi-complete Noetherian local ring $(R,\fm)$ an
$R$-module $D$ is cyclically pure injective if and only if there
is a family $\{C_\lambda\}_{\lambda\in \Lambda}$ of cocyclic
modules such that $D$ is isomorphic to a direct summand of
$\Pi_{\lambda\in \Lambda}C_\lambda$. Also, we prove that over a
local Noetherian ring $(R,\fm)$ every finitely generated
$R$-module $M$ that has small cofinite irreducibles is pure in
every cyclically pure extension of $M$. As a result, we deduce
that over a complete local ring every finitely generated module
which has small cofinite irreducibles is cyclically pure
injective.

\section{Cyclically pure extensions of modules}

Let $\mathcal{S}$ denote the class of all $R$-modules $M$ such
that there are an integer $n\in \mathbb{N}$ and a cyclic submodule
G of $R^n$ such that $M$ is isomorphic to $R^n/G$. In the sequel,
we show that cyclically pure exact sequences and cyclically pure
injective modules are coincide with $\mathcal{S}$-pure exact
sequences and $\mathcal{S}$-pure injective modules, respectively.

\begin{definition} i) Recall that an exact sequence $0\lo A\lo
B\lo C\lo 0$ is said to be {\it cyclically pure} if the natural
map $R/\fa\otimes_RA\lo R/\fa\otimes_RB$ is injective for all
finitely generated ideals $\fa$ of $R$. Also, an $R$-monomorphism
$f:A\lo B$ is said to be {\it cyclically pure}, if the exact
sequence $0\lo A\overset{f}\lo B\overset{nat}\lo B/f(A)\lo 0$ is
cyclically pure. Moreover, a submodule $A$ of an $R$-module $B$ is
called {\it cyclically pure submodule} if the inclusion map
$A\hookrightarrow B$ is cyclically pure.
\\
ii) An $R$-module $D$ is called {\it cyclically pure injective} if
for any cyclically pure exact sequence $0\lo A\lo B\lo C\lo 0$,
the induced sequence $$0\lo \Hom_R(C,D)\lo \Hom_R(B,D)\lo
\Hom_R(A,D)\lo 0$$ is exact.
\end{definition}

In the sequel, we use the abbreviation CP for the term
``cyclically pure".

\begin{proposition} Suppose $0\lo A\hookrightarrow B\overset{\varphi}
\lo C\lo 0$ is an exact sequence of $R$-modules and
$R$-homomorphisms. The following are equivalent:\\
i) For any $M\in \mathcal{S}$, the induced homomorphism
$\Hom_R(M,B)\lo
\Hom_R(M,C)$ is surjective.\\
ii) If the linear equation $\sum_{i=1}^nr_ix_i=a$, $a\in A,
r_1,r_2,\dots, r_n\in R$ is solvable in B, then it is also
solvable in A.\\
iii) $\fa B\cap A=\fa A$ for any finitely generated ideal $\fa$ of $R$.\\
iii)' $\fa B\cap A=\fa A$ for any ideal $\fa$ of $R$.\\
iv) The exact sequence $0\lo A\hookrightarrow
B\overset{\varphi}\lo C\lo 0$ is cyclically pure.\\
iv)' The natural map $R/\fa\otimes_RA\lo R/\fa\otimes_RB$ is
injective for all ideals $\fa$ of $R$.
\end{proposition}

{\bf Proof.}  $ii)\Rightarrow i)$ Let $M=R^n/G$ where
$G=R(r_1,r_2,\dots, r_n)$ is a cyclic submodule of $R^n$ and take
an element $f$ in $\Hom_R(M,C)$. For each $1\leq i\leq n$, set
$m_i=e_i+G$ where $\{e_1,e_2,\dots, e_n\}$ is the standard basis
of $R^n$. There is $y_i\in B$ such that $\varphi(y_i)=f(m_i)$ for
all $1\leq i\leq n$.  One can see easily that
$a:=\sum_{i=1}^nr_iy_i\in A$. Hence, by the assumption there are
$z_1,z_2,\dots z_n\in A$ such that $\sum_{i=1}^nr_iz_i=a$. We
define $g:M\lo B$, by $g(\sum_{i=1}^ns_im_i)=
\sum_{i=1}^ns_i(y_i-z_i)$ for all $s_1,s_2,\dots s_n\in R$.
Suppose that $\sum_{i=1}^ns_im_i=0$ for some $s_1,s_2,\dots s_n\in
R$. Then $(s_1,s_2,\dots s_n)\in G$. Hence there is $b\in R$ such
that $(s_1,s_2,\dots s_n)=b(r_1,r_2,\dots r_n)$, and so
$$\sum_{i=1}^ns_i(y_i-z_i)=b(\sum_{i=1}^nr_iy_i-
\sum_{i=1}^nr_iz_i)=0.$$ Therefore $g$ is well-defined. It is easy
to see that $\varphi g=f$.

$i)\Rightarrow ii)$ Suppose that $n\in \mathbb{N}$, $r_1,r_2,\dots
r_n\in R$ and  $y_1,y_2,\dots y_n\in B$ are such that
$a:=\sum_{i=1}^nr_iy_i$ is an element of $A$. Set
$M=R^n/R(r_1,r_2,\dots r_n)$ and $m_i=e_i+R(r_1,r_2,\dots r_n)$
for all $1\leq i\leq n$. Define $f:M\lo C$, by
$f(\sum_{i=1}^ns_im_i)=\sum_{i=1}^ns_i\varphi(y_i)$. It is a
routine check that $f$ is a well-defined $R$-homomorphism. By the
assumption, there exists an $R$-homomorphism $g:M\lo B$ such that
$\varphi g=f$. We have $\varphi(y_i)=f(m_i)=\varphi(g(m_i))$, and
so $y_i-g(m_i)\in A$ for all $1\leq i\leq n$. Now, we have
\begin{equation*}
\begin{split}
 \sum_{i=1}^nr_i(y_i-g(m_i))&=\sum_{i=1}^nr_iy_i-
\sum_{i=1}^nr_ig(m_i)\\&=a-
g(\sum_{i=1}^nr_im_i)\\&=a-g((r_1,r_2,\dots r_n)+(r_1,r_2,\dots
r_n)R)\\&=a.
\end{split}
\end{equation*}

Next, the equivalences $ii)\Rightarrow iii)'$,  $iii)'\Rightarrow
iv)'$ and $iv)'\Rightarrow iv)$ are clear.  Also, the implications
$iv)\Rightarrow iii)$ and $iii)\Rightarrow ii)$ are obvious, and
so the proof is complete. $\Box$

Let $\{E_i\}_{i\in I}$ be a class of $R$-modules. It is known that
$\prod_{i\in I}E_i$ is an injective $R$-module if and only if
$E_i$ is injective for all $i\in I$. By using the standard
argument of this classical result, we can deduce the following
analogue conclusion for CP-injective modules.

\begin{lemma} Let $\{D_i\}_{i\in I}$ be a class of $R$-modules.
Then $\prod_{i\in I}D_i$ is a CP-injective $R$-module if and only
if $D_i$ is CP-injective for all $i\in I$.
\end{lemma}

\begin{lemma} Let $\fb$ be an ideal of $R$. Then any CP-injective
$R/\fb$-module is also CP-injective as an $R$-module.
\end{lemma}

{\bf Proof.}  Let $D$ be a CP-injective $R/\fb$-module. Assume
that $M$ and $N$ are two $R$-modules. Let $\psi:N\lo M$ be a
CP-homomorphism and let $f\in \Hom_R(N,D)$. Since $\fb D=0$, the
map $f$ induces the $R/\fb$-homomorphism $f^*:N/\fb N\lo D$,
defined by $f^*(x+\fb N)=f(x)$ for all $x+\fb N\in N/\fb N$. From
our assumption on $\psi$, we deduce that the induced
$R/\fb$-homomorphism $\psi^*:N/\fb N\lo M/\fb M$ is cyclically
pure. Thus there is an $R/\fb$-homomorphism $h:M/\fb M\lo D$ with
$h\psi^*=f^*$. Let $g=h\pi$ where $\pi$ is the natural epimorphism
$M\lo M/\fb M$. Then $g\psi=f$. $\Box$

\begin{theorem} Let $M$ be an $R$-module. Then there are a CP-injective
$R$-module $D$ and a CP-homomorphism $\varphi:M\lo D$.
\end{theorem}

{\bf Proof.}  Let $R^*$ denote the set of all finitely generated
ideals of $R$. Set $D=\prod_{\fa \in R^*}E_{R/\fa}(M/\fa M)$ where
$E_{R/\fa}(M/\fa M)$ denotes the injective envelope of the
$R/\fa$-module $M/\fa M$. Define $\varphi:M\lo D$, by
$\varphi(x)=(x+\fa M)_{\fa \in R^*}$ for all $x\in M$. It follows,
by  Lemmas 2.3 and 2.4 that $D$ is CP-injective. Clearly,
$\varphi$ is injective. Next, we prove that $\varphi$ is
cyclically pure. To this end, let $\fb$ be an arbitrary finitely
generated ideal of $R$ and let $y$ be an arbitrary element of $\fb
D\cap \varphi(M)$. Then $y=\varphi(x)$ for some $x\in M$. Since
$\varphi(x)\in \fb D$, it turns out that $x+\fb M\in \fb
E_{R/\fb}(M/\fb M)=0$. Thus $x\in \fb M$, and so $y\in \fb
\varphi(M)$, as required. $\Box$

\begin{corollary} Let $M$ be an $R$-module. There is an extension
$D$ of $M$ such that $D$ is CP-injective and it contains $M$ as a
CP-submodule.
\end{corollary}

{\bf Proof.} Let $M$ be an $R$-module. By Theorem 2.5, there are a
CP-injective $R$-module $D'$ and a CP-homomorphism $\varphi:M\lo
D'$. Using [{\bf 9}, Proposition 1.1], it turns out that there are
an extension $D$ of $M$ and an isomorphism $\psi:D\lo D'$ which is
such that $\psi(x)=\varphi(x)$ for all $x\in M$. It is easy to see
that the inclusion map $M\hookrightarrow D$ is cyclically pure.
$\Box$

\begin{remark} i) The analogue of some of our results for RD-purity were
proved by R.B. Warfield (see e.g. [{\bf 5}, Chapter XIII.1]).\\
ii) One can adapt the method of Warfield's proof of existence of
RD-injective envelopes for proving the existence of CP-injective
envelopes. We present a different proof for existence of
CP-injective envelopes in Section 4.
\end{remark}

\section{Two characterizations}

In this section, we present two characterizations of CP-injective
modules. First, we bring a definition.

\begin{definition} Let $M$ be an $R$-module and $N$ a CP-submodule of $M$.
Then $M$ is called {\it essential CP-extension} of $N$, if there
is not any nonzero submodule $K$ of $M$ such that $K\cap N=0$ and
$(K+N)/K$ is a CP-submodule of $M/K$.
\end{definition}

For a submodule $N$ of an $R$-module $M$, it is known that $M$ is
an essential extension of $N$ if and only if for any $R$-module
$L$ an $R$-homomorphism $\varphi:M\rightarrow L$ is injective
whenever $\varphi\mid_N$ is injective. Similarly, for essential
CP-extensions we have the following characterization.

\begin{lemma} Let $N$ be a CP-submodule of an $R$-module $M$. Then $M$
is essential CP-extension of $N$ if and only if for any
homomorphism $\varphi:M\rightarrow L$ such that $\varphi\mid_N$ is
a CP-homomorphism, it follows that $\varphi$ is injective.
\end{lemma}

{\bf Proof.} Suppose that $M$ is an essential CP-extension of $N$.
Let $\varphi:M\lo L$ be a homomorphism such that $\varphi\mid_N$
is a CP-homomorphism. Let $K=\ker \varphi$ and let
$\varphi^*:M/K\lo L$ denote the natural monomorphism which induced
by $\varphi$. Also, let $\rho:(K+N)/K\lo N$ denote the natural
isomorphism. Note that, because $\varphi\mid_N$ is injective, it
turns out that $K\cap N=0$. From the commutative diagram

\begin{equation*}
\setcounter{MaxMatrixCols}{11}
\begin{matrix}
&(K+N)/K &\hookrightarrow &M/K
\\&\mapdown{\rho}  & &\mapdown{\varphi^*}
& & & &
\\&N &\overset{\varphi\mid_N}\lo
&L &
\end{matrix}
\end{equation*}

we deduce that  $(K+N)/K$ is a CP-submodule of $M/K$. Therefore,
it follows that $K=0$. The proof of the converse is easy and we
leave it to the reader. $\Box$

\begin{lemma} Let $N$ be a CP-submodule of an $R$-module $M$.
Then, there exists a submodule $K$ of $M$ such that \\
i) $K\cap N=0$,\\
ii) $(K+N)/K$ is a CP-submodule of $M/K$, and\\
that $K$ is maximal with respect to inclusion among all submodules
of $M$ which satisfy the conditions i) and ii). In particular,
$M/K$ is an essential CP-extension of $(K+N)/K$.
\end{lemma}

{\bf Proof.} Let $\Sigma$ denote the class of all submodules of
$M$ which satisfy the conditions i) and ii). Then $\Sigma$ is not
empty, because $0\in \Sigma$. Let $\Omega$ be a totally ordered
subclass of $\Sigma$. Set $K=\cup_{K_\alpha\in \Omega}K_\alpha$.
We show that $K$ satisfies the conditions i) and ii). Clearly,
$K\cap N=0$. In view of Proposition 2.2, it is enough to show that
$(K+\fa M)\cap (K+N)\subseteq (\fa N+K)$ for any ideal $\fa$ of
$R$. But, it is a routine check, because by Proposition 2.2
$(K_\alpha+\fa M)\cap (K_\alpha+N)=\fa N+K_\alpha$ for any ideal
$\fa$ of $R$ and all $K_\alpha\in \Omega$. Thus the conclusion
follows by Zorn's Lemma. Now, we prove the last assertion. Assume
there is a submodule $L/K$ of $M/K$ such that $L/K\cap
((N+K)/K)=0$ and that $(L+N)/L$ is a CP-submodule of $M/L$. Then
$$L\cap N\subseteq L\cap (N+K)=K,$$
and so $L\cap N\subseteq K\cap N=0$. Thus $L\in \Sigma$ and so
$L=K$, by the assumption on $K$. Therefore, $M/K$ is an essential
CP-extension of $(K+N)/K$, as required. $\Box$

Now, we present our first characterization of CP-injective
modules.

\begin{theorem}Let $D$ be an $R$-module. Then the following are
equivalent:\\
i) $D$ is CP-injective.\\
ii) For any CP-homomorphism $f:A\rightarrow B$, every homomorphism
from $A$ to $D$ can be extended to a homomorphism
from $B$ to $D$.\\
iii) Every CP-exact sequence $0\lo D\lo M\lo N\lo 0$ splits.\\
iv) $D$ is a direct summand of every $R$-module $L$ which is such
that $D$ is a CP-submodule of $L$.\\
v) $D$ has no proper essential CP-extension.
\end{theorem}

{\bf Proof.} The implications $i)\Rightarrow ii)$, $ii)\Rightarrow
iii)$, and $iii)\Rightarrow iv)$ are clear.

$iv)\Rightarrow v)$ Let $M$ be a essential CP-extension of $D$.
Then, there is a submodule $L$ of $M$ such that $M=L+D$ and $L\cap
D=0$. Since $M$ is essential CP-extension of $D$ and
$(L+D)/L=M/L$, we deduce that $L=0$, and so $M=D$.

$v)\Rightarrow iv)$ Suppose $L$ is a CP-extension of $D$. It may
be assumed that $L$ is a proper CP-extension of $D$. By Lemma 3.3,
there is a submodule $K$ of $L$ such that $L/K$ is an essential
CP-extension of $(D+K)/K$ and that $D\cap K=0$. But $D$ has no
proper essential  CP-extension, so $D+K=L$, from which it follows
that $L=D\oplus K$.

$iv)\Rightarrow i)$ By Corollary 2.6, there exists a CP-injective
extension $L$ of $D$. Therefore, $D$ is a direct summand of $L$,
and so it is CP-injective, by Lemma 2.3. $\Box$

Let $D$ be an $R$-module. In [{\bf 3}, Corollary 2.12], we proved
that $D$ is pure injective if and only if $D$ is isomorphic to a
direct summand of a module of the form $\Hom_R(L,E)$ where $E$ is
an injective $R$-module and $L$ is the direct sum of a family of
finitely generated modules. Next, we will establish a similar
characterization for CP-injective modules. First, we need to the
following lemma.

\begin{lemma} Let $\fa$ be an ideal of $R$. Then an $R/\fa$-module
$D$ is injective as an $R/\fa$-module if and only if there is an
injective $R$-module $E$ such that $D$ is equal to $0:_E\fa$.
\end{lemma}

{\bf Proof.}  The ``if" part is known and it is easy to check. For
the converse, one only need to note that for an injective
$R/\fa$-module $D$, we have $$D=E_{R/\fa}(D)=0:_{E_R(D)}\fa.$$
Here $E_R(D)$ (resp. $E_{R/\fa}(D)$) denotes the injective
envelope of $D$ as an $R$-module (resp. $R/fa$-module). $\Box$

\begin{theorem} Suppose $D$ is an $R$-module. Then the following
are equivalent:\\
i) $D$ is CP-injective.\\
ii) There is a family $\{\fa_\lambda\}_{\lambda\in \Lambda}$ of
finitely generated ideals of $R$ such that $D$ is isomorphic to a
direct summand of an $R$-module of the form $\prod_{\lambda\in
\Lambda}E_\lambda$ where $E_\lambda$ is an injective
$R/\fa_\lambda$-module for all $\lambda\in \Lambda$.\\
iii) $D$ is isomorphic to a direct summand of a module of the form
$\Hom_R(L,E)$ where $E$ is an injective module and $L$ is the
direct sum of a family of finitely presented cyclic modules.
\end{theorem}

{\bf Proof.}  In view of the proof of Theorem 2.5, the equivalence
$i)\Leftrightarrow ii)$ follows by Lemmas 2.3 and 2.4.

$ii)\Rightarrow iii)$ Suppose that there is a family
$\{\fa_\lambda\}_{\lambda\in \Lambda}$ of finitely generated
ideals of $R$ such that $D$ is isomorphic to a direct summand of
an $R$-module of the form $\prod_{\lambda\in \Lambda}E_\lambda$
where $E_\lambda$ is an injective $R/\fa_\lambda$-module for all
$\lambda\in \Lambda$. By Lemma 3.5, for each $\lambda\in \Lambda$,
there is an injective $R$-module $D_\lambda$ such that
$E_\lambda=0:_{D_\lambda}\fa_\lambda$. Hence $\prod_{\lambda\in
\Lambda}E_\lambda=\prod_{\lambda\in
\Lambda}\Hom_R(R/\fa_\lambda,D_\lambda)$. Now, let
$L=\oplus_{\lambda \in \Lambda}R/\fa_\lambda$ and
$E=\prod_{\lambda\in \Lambda}D_\lambda$. Thus $E$ is an injective
$R$-module and $\prod_{\lambda\in \Lambda}E_\lambda$ is a direct
summand of the $R$-module $\Hom_R(L,E)$. Therefore, $D$ is
isomorphic to a direct summand of the $R$-module $\Hom_R(L,E)$.

$iii)\Rightarrow ii)$ It is clear, by Lemma 3.5. $\Box$

\section{CP-injective envelops}

In this section, we show that every R-module possesses a unique,
up to isomorphism, CP-injective envelope.

\begin{definition} i) Let $N$ be an $R$-module. A CP-essential
extension $M$ of $N$ is said to be maximal if there is no proper
extension of $M$ which is a CP-essential extension of $N$.\\
ii) Let $M$ be a CP-submodule of a CP-injective $R$-module $D$. We
say that $D$ is a {\it minimal CP-injective extension} of $M$, if
there is not any proper CP-injective submodule of $D$ containing
$M$.
\end{definition}

\begin{lemma} Let $N$ be an $R$-module and $M$ a CP-essential
extension of $N$. There exists a maximal CP-essential extension
$C$ of $N$ containing $M$.
\end{lemma}

{\bf Proof.} Suppose the contrary is true. By induction on ordinal
numbers, we show that for any ordinal $\beta$, there is a
CP-essential extension $C_\beta$ of $N$ containing $M$. Let
$\beta$ be an ordinal and assume that $C_\alpha$ is defined for
all $\alpha<\beta$. Assume $\beta$ is a predecessor $\beta-1$.
Since $C_{\beta-1}$ is not a maximal CP-essential extension of
$N$, there is a proper extension $C_\beta$ of $C_{\beta-1}$ such
that $C_\beta$ is a CP-essential extension of $N$. If $\beta$ is a
limit ordinal, then it is routine check that
$C_\beta:=\cup_{\alpha<\beta}C_\alpha$ is a CP-essential extension
of $N$. By Corollary 2.6, there is an extension $D$ of $N$ such
that $D$ is CP-injective and it contains $N$ as a CP-submodule.
Let $\beta$ be an ordinal with $\mid \beta \mid>\mid D \mid$.
Then, by Lemma 3.2 the inclusion map $N\hookrightarrow D$ can be
extended to a monomorphism $\psi:C_\beta\lo D$. Hence $\mid \beta
\mid\le\mid C_\beta \mid \leq\mid D \mid$, which is a
contradiction. $\Box$

\begin{lemma} Let $M$ and $M'$ be two $R$-modules and
let $f:M\lo M'$ be an isomorphism. Let $N$ be a submodule $M$ and
$N'=f(N)$.\\
i) $N$ is a CP-submodule of $M$ if and only $N'$ is a CP-submodule
of $M'$.\\
ii) $M$ is a CP-essential extension of $N$ if and only $M'$ is a
CP-essential extension of $N'$.\\
iii) $M$ is a maximal CP-essential extension of $N$ if and only
$M'$ is a maximal CP-essential extension of $N'$.
\end{lemma}

{\bf Proof.} i) is clear.

ii) Assume that $M$ is an CP-essential extension of $N$. By i),
$N'$ is a CP-submodule of $M'$. Let $K'$ be a submodule of $M'$
such that $K'\cap N'=0$ and that $(K'+N')/K'$ is a CP-submodule of
$M'/K'$. Let $K=f^{-1}(K')$. Then $K\cap N=0$, because $f$ is
monic. On the other hand, if $f^*:M/K\lo M'/K'$ denotes the
natural isomorphism induced by $f$, then i) yields that $(K+N)/K$
is a CP-submodule of $M/K$. Note that $f^*((K+N)/K)=(K'+N')/K'$.
Thus $K=0$, and so $K'=0$. Hence $M'$ is a CP-essential extension
of $N'$. The converse follows by the symmetry. Note that
$f^{-1}:M\lo M'$ is an isomorphism with $f^{-1}(N')=N$.

iii) By the symmetry, it is enough to show the ``only if" part.
Suppose that $M$ is a maximal CP-essential extension of $N$. By
ii), $M'$ is a CP-essential extension of $N'$. Let $L'$ be an
extension of $M'$ such that it is a CP-essential extension of
$N'$. By [{\bf 9}, Proposition 1.1], there is an extension $L$ of
$M$ and an isomorphism $g:L\lo L'$ such that the following diagram
commutes.
\begin{equation*}
\setcounter{MaxMatrixCols}{11}
\begin{matrix}
&N &{\hookrightarrow}  &M &{\hookrightarrow} &L
\\&\mapdown{f|_N} & &\mapdown{f}
 & &\mapdown{g} & & & &
\\  &N' &{\hookrightarrow}
&M' &{\hookrightarrow} &L'& &
\end{matrix}
\end{equation*}
It follows by ii), that $L$ is a CP-essential extension of $N$.
Hence $L=M$, by the maximality assumption on $M$. Therefore
$L'=M'$, as required. $\Box$

\begin{corollary} Let $M$ be a CP-injective $R$-module and $N$
a CP-submodule of $M$. There is a submodule $D$ of $M$ which is
maximal CP-essential extension of $N$.
\end{corollary}

{\bf Proof.} By Lemma 4.2, there exists a maximal CP-essential
extension  $L$ of $N$. In view of Lemma 3.2, there is a
monomorphism $\psi:L\lo M$ such that $\psi|_N$ is equal to the
inclusion map $N\hookrightarrow M$. Let $D=\psi(L)$. Since
$\psi:L\lo D$ is an isomorphism, it follows by Lemma 4.3 iii),
that $D$ is also a maximal CP-essential extension of $N$. $\Box$

\begin{proposition} Suppose that $M$ is an $R$-module and that $D$
is a maximal CP-essential extension of $M$. Then $D$ is a
CP-injective $R$-module.
\end{proposition}

{\bf Proof.} In view of Theorem 3.4, it is enough to show that $D$
is a direct summand of every $R$-module which contains $D$ as a
CP-submodule. Let $D$ be a CP-submodule of an $R$-module $L$. By
Lemma 3.3, there exists a submodule $K$ of $L$ such that $K\cap
M=0$ and that $L/K$ is a CP-essential extension of $(K+M)/K$. We
show that $L$ is the direct sum of $K$ and $D$. First, we show
that $K\cap D=0$. Let $K_1=K\cap D$. Then $K_1\cap M=0$ and since
the natural embedding $M\cong (K+M)/K$ in $L/K$ can be factored
throw the natural embedding $M\cong (K_1+M)/K_1$ in $D/K_1$, it
follows that $(K_1+M)/K_1$ is a CP-submodule of $D/K_1$. Thus
$K_1=0$, as required. Now, let $f:D\lo (K+D)/K$ denote that the
natural isomorphism. Then $f(M)=(K+M)/K$. Thus, it follows by
Lemma 4.3 iii) that $(K+D)/K$ is a maximal CP-essential extension
of $(K+M)/K$. But $L/K$ is a CP-essential extension of $(K+M)/K$
and $(K+D)/K\subseteq L/K$. Thus $L=K+D$. $\Box$

Now, we are ready to prove the main result of this section.

\begin{theorem} Let $D$ be an $R$-module and $M$ a submodule of $D$.
The following are equivalent:\\
i) $D$ is a maximal CP-essential extension of $M$.\\
ii) $D$ is a CP-essential extension of $M$ which is CP-injective.\\
iii) $D$ is a minimal CP-injective extension of $M$.
\end{theorem}

{\bf Proof.} $i)\Rightarrow ii)$ is clear by Proposition 4.5.

$ii)\Rightarrow iii)$ Suppose $D'$ is a submodule of $D$
containing $M$ such that $D'$ is CP-injective. By Corollary 4.4,
there exists a submodule $D''$ of $D'$ which is a maximal
CP-essential extension of $M$. Since $D$ is a CP-essential
extension of $M$, it turns out that $D''=D$. Hence $D'=D$.

$iii)\Rightarrow i)$ By Corollary 4.4, there is a submodule $D'$
of $D$ such that $D'$ is a maximal CP-essential extension of $M$.
Now, the module $D'$ is CP-injective, by Proposition 4.5, thus
$D'=D$, by the minimality assumption on $D$. $\Box$

\begin{corollary} Let $M$ be an $R$-module. Then there exists
an $R$-module $D$ satisfying the equivalent conditions i), ii) and
iii) in Theorem 4.6. Moreover, if $D_1, D_2$ are both minimal
CP-injective extensions of $M$ and $i:M\lo D_1$ and $j:M\lo D_2$
denote the related inclusion maps, then there is an isomorphism
$\theta:D_1\lo D_2$ such that $\theta i=j$.
\end{corollary}

{\bf Proof.} The existence of a such $R$-module $D$ follows by
Lemma 4.2. Now, assume that the $R$-modules $D_1$ and $D_2$ are
minimal CP-injective extensions of $M$. In view of Lemma 3.2,
there exists a monomorphism $\theta:D_1\lo D_2$ such that the
following diagram commutes.

\begin{equation*}
\setcounter{MaxMatrixCols}{11}
\begin{matrix}
&M &\overset{i}\hookrightarrow &D_1
\\&\mapdown{id_M}  & &\mapdown{\theta}
& & & &
\\&M &\overset{j}\hookrightarrow
&D_2 &
\end{matrix}
\end{equation*}

The module $\theta(D_1)$ is a CP-injective submodule of $D_2$ that
contains $M$. Hence $\theta(D_1)=D_2$, by the minimality
assumption on $D_1$. This completes the proof. $\Box$

\begin{definition} Let $M$ be an $R$-module and $D$ an extension
of $M$. If $D$ satisfies one of the equivalent conditions of
Theorem 4.6, then $D$ is said  to be the {\it CP-injective
envelope} of $M$.
\end{definition}

Let $\chi$ denote a class of $R$-modules. We recall the notion of
$\chi$-envelope from [{\bf 11}].

\begin{definition} Let $M$ be an $R$-module. An $R$-module
$D\in \chi$ is called $\chi$-envelope of $M$ if there is a
homomorphism
$\varphi:M\lo D$ such that\\
i) for any homomorphism $\varphi^{\prime}:M\lo D^{\prime}$ with
$D^{\prime}\in \chi$, there is a homomorphism $f:D\lo D^{\prime}$
such with $\varphi^{\prime}=f\varphi$, and\\
ii) if a homomorphism $f:D\lo D$ is such that $\varphi=f\varphi$,
then $f$ must be an automorphism.
\end{definition}

By [{\bf 11}, Proposition 1.2.1], the $\chi$-envelope of an
$R$-module is unique up to isomorphism. Now, we present our last
result.

\begin{theorem} Let $\chi$ be the class of all CP-injective $R$-modules
 and $M$ an $R$-module. Let $D$ be a CP-injective envelope of $M$.
Then $D$ is isomorphic to the $\chi$-envelope of $M$.
\end{theorem}

{\bf Proof.} Let $\varphi:M\lo D$ denote the inclusion
homomorphism. Let $D^{\prime}\in \chi$ and $\psi:M\lo D^{\prime}$
be a homomorphism. By Theorem 3.4, there exists a homomorphism
$f:D\lo D^{\prime}$ such that $\psi=f\varphi$.

Now, suppose a homomorphism $f:D\lo D$ is such that
$\varphi=f\varphi$. By Lemma 3.2, the map $f$  is injective. By
Lemma 4.3 iii), $f(D)$ is also a maximal CP-essential extension of
$M$. Hence $f(D)$ is CP-injective, by Proposition 4.5. On the
other hand, we have $M\subseteq f(D)\subseteq D$. Therefore, by
using Theorem 4.6, we deduce that $f(D)=D$, and so $f$ is an
automorphism, as required. $\Box$

\section{Cocyclic modules}

In [{\bf 3}], we showed that over a Noetherian ring $R$ an
$R$-module $D$ is pure injective if and only if $D$ is isomorphic
to a direct summand of the direct product of a family of Artinian
modules. In this section, we intent to provide an analogue
characterization for CP-injective modules, by using cocyclic
modules instead of Artinian modules. It is known that any Artinian
module is pure injective, but it is not the case that every
cocyclic module is CP-injective (see Example 5.6). Thus, it is
interesting to know when a cocyclic modules is CP-injective.
First, we recall some definitions.

\begin{definition} i) (See [{\bf 4}, page 4]) An $R$-module $M$ is
called {\it cocyclic} if $M$ is isomorphic to a submodule of the
injective envelope of a simple module.\\
ii) (See [{\bf 1}]) An $R$-module $M$ is called {\it subdirectly
irreducible} if for any family $\{M_i\}_{i\in I}$ of $R$-modules
and any monomorphism $f:M\lo \Pi_{i\in I}M_i$, there exists $i\in
I$ such that the map $\pi_if:M\lo M_i$ is injective, where
$\pi_i:\Pi_{i\in I}M_i\lo M_i$ denotes the $i$-th projection map.
\end{definition}

In the following result, we collect some other conditions that are
equivalent to the definition of a cocyclic module.

\begin{proposition} Let $M$ be a nonzero $R$-module. Then the following
are equivalent:\\
i) $M$ is cocyclic.\\
ii) $E_R(M)=E_R(S)$ where $S$ is a simple module.\\
iii) The socle of $M$ is simple and $M$ is an essential extension
of its socle.\\
iv) The intersection of all nonzero submodules of $M$ is nonzero.\\
v) There exists an element $c\in M$ such that for every $R$-module
$N$ and every $R$-homomorphism $f:M\lo N$, it follows that $f$ is
injective if and only if $c\not\in \ker f$.\\
vi) The intersection of all nonzero submodules of $M$ is a simple
 submodule of $M$.\\
vii) $M$ is subdirectly irreducible.
\end{proposition}

{\bf Proof.} The equivalence $i)\Leftrightarrow ii)$ follows by
[{\bf 9}, Proposition 2.28]. Also, the equivalences
$iv)\Leftrightarrow v)$ and $iv)\Leftrightarrow vi)$ are clear.

$i)\Rightarrow iii)$ Suppose the simple $R$-module $S$ is such
that $M$ is isomorphic to a submodule of $E_R(S)$. Then $M$
possesses a simple submodule $S'$ such that every nonzero
submodule of $M$ contains $S'$. Hence, the socle of $M$ is simple.
On the other hand, $M$ is essential extension of its socle, by
[{\bf 9}, Proposition 3.17].

Next, the equivalence $iii)\Leftrightarrow vi)$ and the
implication $iii)\Rightarrow ii)$ both are deduced, by [{\bf 9},
Proposition 3.17].

$vi)\Rightarrow vii)$ Consider the family $\{M_i\}_{i\in I}$ of
$R$-modules and a monomorphism $f:M\lo \Pi_{i\in I}M_i$. For each
$i\in I$, let $\pi_i:\Pi_{i\in I}M_i\lo M_i$ denote the $i$-th
projection map. Assume that the simple $R$-module $S$ is equal to
the intersection of all nonzero submodules of $M$ and let $x$ be a
nonzero element of $S$. Since $f(x)\neq 0$, it follows that there
is $i\in I$ such that $(\pi_if)(x)\neq 0$. This implies that
$\ker(\pi_if)=0$, because otherwise $S\subseteq \ker(\pi_if)$,
which a contradiction.

Finally, we prove that vii) implies iv). Let
$\{N_{\lambda}\}_{\lambda\in \Lambda}$ denote the set of all
nonzero submodules of $M$ and let $f:M\lo \Pi_{\lambda\in
\Lambda}M/N_{\lambda}$ denote the natural homomorphism defined by
$x\mapsto (x+N_\lambda)_{\lambda\in \Lambda}$. Denote
$\bigcap_{\lambda\in \Lambda}N_\lambda$ by $S$. If $S=0$, then $f$
is injective, and so there is  $\lambda\in \Lambda$ such that
$\pi_{\lambda}f:M\lo M/N_{\lambda}$ is injective. This implies
that $N_\lambda=0$, which is a contradiction. $\Box$

\begin{proposition} Let $M$ be an $R$-module and let
 $\{N_i\}_{i\in I}$ denote the set of all submodules $N$ of $M$,
  such that $M/N$ is cocyclic. Then the natural map $\psi:M\lo \Pi_{i
 \in I}M/N_i$ is cyclically pure. In particular, if $M$ is CP-injective
 then $M$ is isomorphic to a direct summand of $\Pi_{i\in I}M/N_i$.
\end{proposition}

{\bf Proof.}  Let $L=\Pi_{i\in I}M/N_i$ and for each $i\in I$ let
$\pi_i:L\lo M/N_i$ denote the $i$-th natural projection map.
Define $\psi:M\lo L$ by $x\mapsto (\pi_i(x))_i$. We show that
$\psi$ is a CP-homomorphism. To this end, let $\fa$ be an ideal of
$R$ and consider the following commutative diagram in which all
maps are natural ones.

\begin{equation*}
\setcounter{MaxMatrixCols}{11}
\begin{matrix}
&M\otimes_RR/\fa &\overset{\psi\otimes id_{R/\fa}}\rightarrow
&L\otimes_RR/\fa
\\&\mapdown{\cong}  & &\mapdown{}
& & & &
\\&M/\fa M &\overset{\theta}\rightarrow
&\Pi_{i\in I}M/(\fa M+N_i) &
\end{matrix}
\end{equation*}

It suffices to show that the bottom map is injective. Let
$\alpha=x+\fa M$ be a nonzero element of $M/\fa M$. Using Zorn's
Lemma, we deduce that  there is a submodule $N$ of $M$ such that
$\fa M\subseteq N$ and $x\notin N$, but $x$ belongs to any
submodule of $M$ which strictly contains $N$. Now, by Proposition
5.2, it turns out that $M/N$ is cocyclic. So, there is $j\in I$
such that $N=N_j$. Since $\pi_j(x)\neq 0$, it follows that
$\theta$ is monomorphism, as required. $\Box$

\begin{definition}(See e.g. [{\bf 7}]) A Noetherian local ring $(R,\fm)$
is called {\it quasi-complete} if for any decreasing sequence
$\{\fa_i\}_{i\in I}$ of ideals of $R$ and any $n\geq 0$, there
exists ${i\in I}$ such that $\fa_i\subseteq (\cap_{i\in
I}\fa_i)+\fm^n$.
\end{definition}

Now, we are ready to present our last characterization of
CP-injective modules.

\begin{theorem} Let $(R,\fm)$ be a quasi-complete local ring.
An $R$-module $D$ is CP-injective if and only if $D$ is isomorphic
to a direct summand of the direct product of a family of cocyclic
modules.
\end{theorem}

{\bf Proof.} Let $E=E_R(R/\fm)$. By [{\bf 7}, Remark 3.2], every
cocyclic $R$-module has the form $0:_E\fa$ for some ideal $\fa$ of
$R$. Thus, by Lemma 2.4 every cocyclic $R$-module is CP-injective.
Now, the conclusion follows by Lemma 2.3 and Proposition 5.3.
$\Box$

\begin{example} By [{\bf 4}, Theorem 6], a Pr\"{u}fer domain $R$ is
locally almost maximal if and only if every cocyclic $R$-module is
pure injective. On the other hand, by [{\bf 2}, Example 2.4] there
exists a valuation domain $R$ such that $R$ is not almost maximal.
Hence cocyclic modules are not CP-injective in general, and so the
converse of Proposition 5.3 is not true.
\end{example}

In [{\bf 6}], M. Hochster, investigated the structure of
Noetherian rings $R$ with the property that $R$ is pure in each
CP-extension of $R$. Let $(R,\fm)$ be a Noetherian local ring, he
defined a finitely generated $R$-module $M$ to have {\it small
cofinite irreducibles} if for every $n\in\mathbb{N}$ there is an
irreducible submodule $Q$ of $M$ such that $Q\subseteq \fm^n M$
and $M/Q$ is Artinian. He showed that a Noetherian ring $R$ is
pure in each CP-extension of $R$ if and only if $R_{\fm}$ has
small cofinite irreducibles for all maximal ideals $\fm$ of $R$.
In this section, we will prove that over a local Noetherian ring
$R$ every finitely generated $R$-module $M$ that has small
cofinite irreducibles is pure in every CP-extension of $M$. As a
result, we deduce that over a complete local ring every finitely
generated module which has small cofinite irreducibles is
CP-injective.

\begin{lemma} Let $R$ be a Noetherian ring and $D$ a finitely generated
cocyclic $R$-module. Then $D$ is CP-injective.
\end{lemma}

{\bf Proof.} There is a maximal ideal $\fm$ of $R$ such that $D$
is isomorphic to a submodule of $E:=E_R(R/\fm)$. Then it is easy
to see that the natural map $D\lo D_{\fm}$ is an isomorphism.
Also, one can check easily that, if D is CP-injective as an
$R_{\fm}$-module then it is also CP-injective as an $R$-module.
So, we may and do assume that $R$ is local with the maximal ideal
$\fm$.

Let $\hat{R}$ denote the completion of $R$ with respect to
$\fm$-adic topology. Each element of $E_R(R/\fm)$ is annihilated
by some power of $\fm$. Hence $E_R(R/\fm)$ has a natural structure
as an $\hat{R}$-module. Note that, if we regard this
$\hat{R}$-module as an $R$-module by means of the natural ring
homomorphism $R\lo \hat{R}$, then we recover the original
$R$-module structure on $E_R(R/\fm)$. Note also that a subset of
$E_R(R/\fm)$ is an $R$-submodule if and only if it is an
$\hat{R}$-submodule. Set $\fa:=\Ann_RD$ and $E:=E_R(R/\fm)$.
Since, D is finitely generated, it turns out that
$\Ann_{\hat{R}}D=\fa \hat{R}$. Therefore, by [{\bf 9}, page 154,
Corollary] we have
$$D=(0:_E\Ann_{\hat{R}}D)=(0:_E\fa).$$
Thus the claim follows, by Lemma 2.4. $\Box$

\begin{theorem} Let $(R,\fm)$ be a Noetherian local ring and $M$ a
finitely generated $R$-module. If $M$ has small cofinite
irreducibles, then $M$ is pure in every CP-extension of $M$.
\end{theorem}

{\bf Proof.} Let $\{N_i\}_{i\in I}$ denote the set of all
submodules $N$ of $M$, such that $M/N$ is cocyclic. Let $L$ and
$\psi:M\lo L$ be as in the proof of Proposition 5.3. Let $C$ be a
CP-extension of $M$ and let $i:M\hookrightarrow C$ denote the
inclusion map. From Lemma 5.7, it follows that $L$ is
CP-injective, and so there is a homomorphism $f:C\lo L$ such that
$fi=\psi$. Therefore, to prove $M$ is pure in $C$, it suffices to
show that $\psi$ is pure. So, we are going to show that for any
finitely generated $R$-module $N$, the induced map
$$\psi\otimes id_N:M\otimes_RN\lo L\otimes_RN$$ is
injective. Assume there exists $n\in\mathbb{N}$ such that $\fm^n
N=0$. Then, since $M$ has small cofinite irreducibles, there is an
irreducible submodule $Q_0$ of $M$ such that $Q_0\subseteq \fm^n
M$ and $M/Q_0$ is Artinian. Then there is $j\in I$ such that
$Q_0=N_j$. For each $i\in I$, let $\pi_i:M\lo M/N_i$ denote the
natural epimorphism. Now, because the modules $M\otimes_RN$ and
$M/N_j\otimes_RN$ are naturally isomorphic, it turns out that
$\pi_j\otimes id_N$ is an isomorphism. Consider the following
commutative diagram.

\begin{equation*}
\setcounter{MaxMatrixCols}{11}
\begin{matrix}
&M\otimes_RN &\overset{\psi\otimes id_N}\rightarrow &L\otimes_RN
\\&\mapdown{id_{M\otimes_RN}}  & &\mapdown{\cong}
& & & &
\\&M\otimes_RN &\overset{\Pi(\pi_i\otimes id_N)}\rightarrow
&\Pi_{i\in I}(M/N_i\otimes_RN) &
\end{matrix}
\end{equation*}
Hence $\psi\otimes id_N$ is injective.

Next, assume that $N$ is an arbitrary finitely generated
$R$-module. Suppose that $\ker(\psi\otimes id_N)$ contains a
nonzero element $x$. Set $K=M\otimes_RN$. Since
$\cap_{i\in\mathbb{N}}\fm^i K=0$, it follows that there is $n\in
\mathbb{N}$ such that $x\notin \fm^n K$. Set $\overline{N}=N/\fm^n
N$ and let $\pi:N\lo \overline{N}$ denote the natural epimorphism.
Because the modules $K/\fm^n K$ and $M\otimes_R\overline{N}$ are
naturally isomorphic, it turns out that the element
$(id_M\otimes\pi)(x)$ of the module $M\otimes_R\overline{N}$ is
nonzero. From the commutative diagram

\begin{equation*}
\setcounter{MaxMatrixCols}{11}
\begin{matrix}
&M\otimes_RN &\overset{id_M\otimes \pi}\rightarrow
&M\otimes_R\overline{N}
\\&\mapdown{\psi\otimes id_N}  & &\mapdown{\psi\otimes
id_{\overline{N}}} & & & &
\\&L\otimes_RN &\overset{id_{L}\otimes \pi}\rightarrow
&L\otimes_R\overline{N} &
\end{matrix}
\end{equation*}
we deduce that $\psi\otimes id_{\overline{N}}$ is not injective,
which is a contradiction in view of the first paragraph of the
proof. $\Box$

\begin{corollary} Let $(R,\fm)$ be a Noetherian complete local ring and
$M$ a finitely generated $R$-module. If $M$ has small cofinite
irreducibles, then $M$ is CP-injective.
\end{corollary}

{\bf Proof.} By Proposition 5.8, $M$ is pure in every CP-extension
of $M$. Thus by [{\bf 7}, Theorem 3.3], $M$ is CP-injective.
$\Box$

\begin{remark} i) Let $R$ be a field. Clearly, every monomorphism is
split and so it is pure. We show that $M=R\oplus R$ doesn't have
small cofinite irreducibles. Suppose the contrary is true. Then
there is an irreducible submodule $Q$ of $M$ such that $Q\subseteq
0M=0$. That is the zero submodule of $M$ is irreducible. Therefore
we achieved at a contradiction. This shows that the converse of
Theorem 5.8 and Corollary 5.9 don't hold. Thus one may consider
these results as generalizations of [{\bf 6}] and  [{\bf 7},
Corollary 3.4],
respectively.\\
ii) It might be interesting to know when the converse of the last
part of Proposition 5.3 holds. Clearly, this is the case when
every cocyclic $R$-module is CP-injective. By Theorem 3.4 and
Proposition 5.2, it is easy to see that if every cocyclic
$R$-module is CP-injective, then the only CP-submodules of a
cocyclic $R$-module are the trivial ones. Is the converse true?
\end{remark}

{\bf Acknowledgments.} We thank the referee for very careful
reading of the manuscript and also for his/her useful suggestions.

%%%%%%%%%%%%%%%%%%%%%%%%%%%%%%%%%%%%%%%%%%%%%%%%%%%%%%%%%%%%%%%%%%%%%%%%%%%%%


\begin{thebibliography}{99}
\bibitem{} F.W. Anderson and K.R. Fuller, {\it Rings and categories of
modules}, 2nd edition, Springer-Verlag, New York, 1992.
\bibitem{} W. Brandal, {\it Almost maximal integral domains and finitely
generated modules}, Trans. Amer. Math. Soc., {\bf 183} (1973),
203-222.
\bibitem{} K. Divaani-Aazar, M.A. Esmkhani and M. Tousi, {\it Two
characterzations of pure injective modules}, Proc. Amer. Math.
Soc., to appear.
\bibitem{} L. Fuchs and A. Meijer, {\it Note on modules
over Pr\"{u}fer domains}, Math. Pannon, {\bf 2}(1) (1991), 3-11.
\bibitem{} L. Fuchs and L. Salce, {\it Modules over non-Noetherian domains},
Mathematical Surveys and Monographs, {\bf 84}, Amer. Math. Soc.,
Providence, RI, 2001.
\bibitem{} M. Hochster, {\it Cyclic purity versus purity in excellent
Noetherian rings}, Trans. Amer. Math. Soc., {\bf 231}(2) (1977),
463-488.
\bibitem{} L. Melkersson, {\it Small cofinite irreducibles
}, J. Algebra, {\bf 196}(2) (1997), 630-645.
\bibitem{} M. Prest, {\it Model theory and modules}, London Mathematical
Society Lecture Note Series, {\bf 130}, Cambridge University
Press, Cambridge, 1988.
\bibitem{} D.W. Sharpe and P. V$\acute{a}$mos, {\it Injective modules},
Cambridge Tracts in Mathematics and Mathematical Physics, {\bf
62}, Cambridge University Press, London-New York, 1972.
\bibitem{} R.B. Warfield, {\it Purity and algebraic compactness for modules},
Pacific J. Math., {\bf 28}, (1969), 699-719.
\bibitem{} J. Xu, {\it Flat covers of modules}, Lecture Notes in Mathematics,
{\bf 1634}, Springer-Verlag, Berlin, 1996.

\end{thebibliography}
\end{document}